\documentclass[a4paper,twoside,final]{article}
\usepackage[latin1]{inputenc}
\usepackage[T1]{fontenc}
\usepackage{amssymb}
\usepackage{dochier}
\usepackage{lpretex}
\usepackage{title}
\usepackage{a4}
\usepackage{amscd}
\EndOfDump
\usepackage{url}

\def\DHpartformat#1{(#1) }
\def\DHrefpart#1{(\DHRefpart{#1})}
\LBsetstyleof{partno}{arabic}

\renewcommand\labelenumi{(\alph{enumi})}

\let\define\def

  \def\F {{\mathbb F}}

\def\N {{\mathbb N}}  \def\P {{\mathbb P}} 
\def\Q {{\mathbb Q}} 
 
 \def\W {{\mathbb W}} 
\def\Z {{\mathbb Z}} 

\define \n {\mathbb N}
\define \z {\mathbb Z}
\define \q {\mathbb Q}
\define \PP {\mathbb P}

\def\sA {{\Cal A}}  
 \def\sE {{\Cal E}} \def\sF {{\Cal F}}
\def\sG {{\Cal G}}  
  
  \def\sO {{\Cal O}}

  \def\sX {{\Cal X}}

\define \cN {\Cal N}
\define \cf {\Cal F}
\define \cg {\Cal G}
\define \cE {\Cal E}
\define \ce {\Cal E}
\define \cc {\Cal C}
\define \cV {\Cal V}
\define \cA {\Cal A}
\define \cK {\Cal K}
\define \cO {\Cal O}
\define \cF {\Cal F}
\define \cn {\Cal N}
\define \cI {\Cal I}
\define \sP {\Cal P}

\define \x {\xi}
\define \y {\eta}
\define \G {\Gamma}
\define \r {\rho}
\define \w {\omega}

\define \eps {\epsilon}

\define \tH {\widetilde H}
\define \tG {\widetilde{\Gamma}}
\define \tW {\widetilde W}
\define \tF {\widetilde F}
\define \tm {\tilde m}
\define \St {\widetilde S}
\define \Xt {\widetilde X}
\define \tS {\widetilde S}
\define \tpsi {\tilde \psi}
\define \tL {\widetilde L}
\define \tE {\widetilde E}
\define \tl {\tilde l}
\define \tA {\widetilde A}
\define \tom {\tilde\omega}
\define \tT {\widetilde T}
\define \tB {\widetilde B}
\define \tf {\tilde f}
\define \tsA {\widetilde{\sA}}
\define \tsF {\widetilde{\sF}}

\define \tM {\widetilde M}
\define \tphi {\widetilde{\phi}}
\define \tR {\widetilde R}
\define \tQ {\widetilde Q}

\def \hR{\widehat R}
\def \hot{\widehat\otimes}

\def \hGa{\widehat\Gamma}
\def \hG{\widehat G}

\def \hS{\widehat S}
\def \hX{\widehat X}

\def\la {\langle} 
\def\ra {\rangle} 
\def\pd {\partial}

\def \Dx1 {\frac{\pd}{{\pd} x_1}}
\def \Dy1 {\frac{\pd}{{\pd} y_1}}
\def \Dz1 {\frac{\pd}{{\pd} z_1}}
\def \Dx2 {\frac{\pd}{{\pd} x_2}}
\def \Dy2 {\frac{\pd}{{\pd} y_2}}
\def \Dz2 {\frac{\pd}{{\pd} z_2}}

\def\q {\quad}

\def\mapdiagr#1{\Big\searrow\rlap{$\raise 5pt\vbox{{\hbox{$\mkern -15mu\scriptstyle#1$}}}$}}   

\def\mapdiagl#1{\llap{$\raise 5pt\vbox{{\hbox{$\scriptstyle#1\mkern
-15mu$}}}$}\Big\swarrow}              

\def\Mapdiagr#1{\nearrow\rlap{$\lower 5pt\vbox{{\hbox{$\mkern
-15mu\scriptstyle#1$}}}$}} 

\def\Mapdiagl#1{\llap{$\lower 5pt\vbox{{\hbox{$\scriptstyle#1\mkern
-15mu$}}}$}\searrow} 

\def\Mapswr#1{\swarrow\rlap{$\lower 5pt\vbox{{\hbox{$\mkern
-15mu\scriptstyle#1$}}}$}}              

\def\Mapnwl#1{\nwarrow\rlap{$\lower 5pt\vbox{{\hbox{$\mkern
-15mu\scriptstyle#1$}}}$}}

\define \Rhook {\hookrightarrow}

\def \half {\raise1pt\hbox{$\scriptstyle
        \frac{1}{2}\displaystyle$}}

\def \x{{\sl X}\llap{$\mkern -2mu {\scriptstyle -}$}}
\def \Hom {\operatorname{Hom}}

\let\Spec\Sp

\define \Kod {\operatorname{Kod}}
\define \dimension {\operatorname{dim}}
\define \codim {\operatorname{codim}}
\define \contr {\operatorname{contr}}
\define \rk {\operatorname{rank}}
\define \im {\operatorname{im}}
\define \Mor {\operatorname{Mor}}
\define \Cl {\operatorname{Cl}}
\define \Hilb {\operatorname{Hilb}}
\define \degree {\operatorname{deg}}
\define \mult {\operatorname{mult}}
\define \Aut {\operatorname{Aut}}
\define \NS {\operatorname{NS}}
\define \Gal {\operatorname{Gal}}
\define \ch {\operatorname{char}}
\define \Jac {\operatorname{Jac}}
\define \Km {\operatorname{Km}}
\define \Sec {\operatorname{Sec}}
\define \Stab {\operatorname{Stab}}
\define \Br {\operatorname{Br}}
\define \inv {\operatorname{inv}}
\define \tr {\operatorname{tr}}
\define \Frob {\operatorname{Frob}}
\define \Symn {\operatorname{Sym}^n}
\define \Ev {\sE^\vee}
\define \ordp {\operatorname{ord}_p}
\define \Supp {\operatorname{Supp}}
\define \Ann {\operatorname{Ann}}
\define \disc {\operatorname{disc}}
\define \Lie {\operatorname{Lie}}
\define \embdim {\operatorname{embdim}}

\def \Der {\operatorname{Der}}

\def\hCL{{\widehat C}_\Lambda}

\def\hot{\hat{\otimes}}

\def\maxid{{\frak{m}}}

\newcommand\piF{\overline{F}}
\newcommand\piG{\overline{G}}
\newcommand\piH{\overline{H}}

\define\cC {\Cal C}

\def\map#1#2#3{\ensuremath{{#1}\co{#2}\to{#3}}}
\def\DP#1{{[#1]}}

\def\W{\cansymb{W}}

\def\hod#1#2#3#4{\ensuremath{\if#30 H^{#2}({#1},{\cal O}_{#1}) \else 
 H^{#2}(#1,\Omega^{#3}\if\relax{#4}\relax_{#1}\else _{#1/#4}\fi)\fi}}

\newcommand\Witt{\W}
\newcommand\Cohen{\Witt}
\newcommand\Reg{\symb{Reg}}

\DocVersion[Nick's original]{1}
\DocVersion[Divided powers and Hirokado's threefold;T.E.;Mon May 14 2001]{2}
\DocVersion[Schroers threefolds;T.E.;Mon May 26 2003]{3}
\DocVersion[Pure deformation theory;T.E.;Wed Jun 18 2003]{4}

\begin{document}
\title{Tangent lifting of deformations in mixed characteristic}

\author{T. Ekedahl}
\address{Department of Mathematics\\
 Stockholm University\\
 SE-106 91  Stockholm\\
Sweden}
\email{teke@matematik.su.se}
\author{N. I. Shepherd-Barron}
\address{D.P.M.M.S.\\
 16 Mill Lane\\
 Cambridge CB2 1SB\\
 U.K.}
\email{nisb@pmms.cam.ac.uk}


\begin{abstract}
This article presents a new approach to the unobstructedness result for
deformations of Calabi-Yau varieties by introducing the \emph{tangent lifting
property} for a functor on artinian local algebras. The verification that the
deformation functor of a Calabi-Yau variety in characteristic $0$ fulfills the
tangent lifting property uses, just as the verification of the $T^1$-lifting
property, the degeneration of the Hodge to de Rham spectral sequence. Our main
use of our methods is however to the mixed characteristic case. In that case to
be able to verify the conditions needed one needs an extension of the criterion
introducing divided powers.
\end{abstract}
\maketitle
Fix a complete noetherian local ring $\Lambda$ with residue field $\k$, and
denote by $C_\Lambda$, respectively $\hCL$, the category of Artin, respectively
complete, local $\Lambda$-algebras with residue field $\k$.  Suppose that $F$ is
a semi-homogeneous cofibered groupoid (cf.\ \cite{SGA7I}) over $C_\Lambda$ such
that $F(\k)$ is equivalent to a point so that $\piF:=\pi_0F$, the functor of
isomorphism classes of objects of $F$, has a hull $R$, say. Define
$A_n=\k[t]/(t^{n+1})$ and $B_n=A_n[\epsilon]/(\epsilon^2)$.  Then $F$ is said to
have the \emph{$T^1$-lifting property} if the natural map $F(B_{n+1})\to
F(B_n)\times_{F(A_n)}F(A_{n+1})$, where the fibre product is the $2$-fibre
product, is essentially surjective.  Ran (\cite{ran92::defor}), Kawamata
(\cite{kawamata92::uenob}) and Fantechi and Manetti (\cite{fantmane99::oen+t})
have shown that if $\Lambda =\k$ and $\ch \k=0$, then $R$ is formally smooth
over $\k$ if $F$ has the $T^1$-lifting property. (The $T^1$-lifting property can
also be phrased as follows. Given $X\in F(A)$, define $T^1(X/A)$ as the set of
isomorphism classes of pairs $(Y,\psi)$, where $Y\in F(A[\eps])$ and $\psi$ is
an isomorphism \map{\psi}{Y_{|A}}{X}.  Then the $T^1$-lifting property is
equivalent to the surjectivity of $T^1(X_n/A_n)\to T^1(X_{n-1}/A_{n-1})$.)
Moreover, the $T^1$-lifting property is satisfied by, e.g., deformations of
Calabi-Yau varieties.

In this article we propose an approach to the question of whether $R$ is
formally smooth that puts different conditions on the functor. Even though this
approach gives a new proof of the unobstructedness for deformations of
Calabi-Yau varieties we shall mainly be interested in the case when $\Lambda$
has mixed characteristic or has positive characteristic.  Roughly speaking, in
the mixed case this involves proving smoothness by (1) lifting the whole tangent
space at once, rather than extending a given tangent direction, and (2) assuming
the existence of some lifting to char. zero.  For deformations of Calabi-Yau
varieties (defined as smooth varieties with trivial dualizing sheaf) in
characteristic zero the fact that these conditions are fulfilled follows in the
same way as the $T^1$-lifting property but make the verification of smoothness
simpler. They are inspired by the fact that in characteristic zero lifting of
tangent vectors to vector fields implies smoothness, a fact that is central in,
for instance, Cartier's proof that a group scheme in characteristic zero is
smooth. Our approach has the added advantage of not requiring that $\piF$ comes
from a cofibered groupoid (in practice this seems almost always to be the case
so it is not clear how great an advantage that is) and hence unless explicitly
stated we shall assume only that $\piF$ is a functor from $C_\Lambda$ to the
category of sets. We shall say that $\piF$ has \emph{the tangent lifting
property} (abbreviated to TLP) if for all $A\in C_\Lambda$ the natural map
$\piF(A[\epsilon])\to \piF(A)\times \piF(\k[\epsilon])$ is surjective.
\begin{remark}
\part It is not clear to us what the precise logical relations between the
notions of $T^1$-lifting and tangent lifting are, although in characteristic zero
$T^1$-lifting implies tangent lifting {\emph{a posteriori}} as tangent lifting
is true when the hull is smooth. It should also be noted that the reason that we
do not need to assume that $\piF$ has an underlying groupoid for the tangent
lifting property is that contrary to the case of $T^1$-lifting only the product
$\piF(A)\times \piF(\k[\epsilon])$ and no fibre product is involved.

\part Schr\"oer \cite{schroeer01::t1} has results similar to the ones presented
here, but with stronger assumptions. Just as our approach gives an alternative in
the characteristic zero to previous ones our method gives a new approach independent of
Schr\"oer's in the mixed characteristic case.  See the remark after Theorem C
below for a more thorough explanation.
\end{remark}
A first re-formulation of this condition is given in the following lemma.
\begin{lemma}\label{T1}
When $\piF$ is $\pi_0F$ of a cofibered groupoid the tangent lifting property is
equivalent to the surjectivity of $T^1(X/A)\to T^1(X_0/\k)$ for all $X\in F(A)$.
\begin{proof}
Consider the restriction map $\pi:F(A[\eps])\to F(A)$. Then for any $X\in
F(A)$, there is a forgetful map $T^1(X/A)\to \pi^{-1}(X)$ and this map is
surjective. Since $T^1(X_0/\k)$ is identified with $F(\k[\eps])$, the lemma is now
(even more) obvious.
\end{proof}
\end{lemma}
If $X_0$ is a smooth variety either in char.\ zero or in char.\ $p$ and $F$
is the deformation groupoid we may identify $T^1(X/A)$ with
$H^1(X,T_{X/A})$, which is isomorphic to $H^1(X,\Omega^{n-1}_{X/A})$
provided that $\omega_{X/A}$ is trivial. Thus the tangent lifting condition may be
analysed through Hodge cohomology. In particular it is always fulfilled for 
Calabi-Yau varieties (that is, those for which $\omega$ is trivial) in
characteristic zero. More generally, the deformation functor of a complete
smooth variety $X_0$ has the tangent lifting property if and only if for every
deformation \map{f}{X}{\Spec A} the specialization homomorphism
$R^1f_*T_{X/A}\to H^1(X_0,T_{X_0})$ is surjective.

\begin{notation} 
We shall abuse language by failing to distinguish between smoothness and formal
smoothness, derivations and continuous derivations, and differentials and
continuous differentials. We shall furthermore, for a sheaf of rings $\sO$ let
$\sO[\epsilon]$ denote $\sO[x]/(x^2)$, where $\epsilon$ is the residue of
$x$. Similarly if $X$ is a scheme we put $X[\epsilon]:=\Sp\sO[\epsilon]$.

To avoid confusion in the case where the involved schemes are non-reduced let us
make explicit that we say that a map \map fXS of (formal) schemes is
\Definition{dominant} if the induced map $f^{-1}\sO_S \to \sO_X$ is injective.
\end{notation}
We shall show that the tangent lifting property for $\piF$ is equivalent to a
tangent lifting property for a deformation hull, namely that any tangent vector
(at the closed point) lifts to a vector field. In characteristic zero this
immediately implies that $R$ is formally smooth whereas in positive
characteristic it implies that defining equations for $R$ can be assumed to be
$p$'th powers. Just as in \cite{schroeer01::t1}, in the case of mixed
characteristic we shall need to have at least one lifting to characteristic
zero. Since this involves rings in $\hCL$ rather than $C_\Lambda$, we make a
convention: Given a functor $\piF:C_\Lambda\to{\underbar{Sets}}$ and $R\in\hCL$,
we write $R=\ili A_\alpha$, with $A_\alpha\in C_\Lambda$, and then define
$\piF(R)$ to be the set of pro-objects $\{\xi_\alpha\in\piF(A_\alpha)\}$.


\begin{Theorem}[A]\label{Main Theorem} Assume that $\piF$ has
the tangent lifting property and let $R$ be a hull for it.

\part \label{char 0} If $\Lambda$ contains $\Q$, then $R$ is smooth.

\part \label{char p} If $\Lambda = \k$ a field and $\ch \k =p>0$, then $R\cong
\pow[\k]{x_1,\dots,x_n}/(f_1,\dots,f_r)$ for some $f_i$ in the maximal ideal of
$\pow[\k]{x_1,\dots,x_n}$ that are of the form $f_i=g_i(x_1^p,\dots,x_n^p)$. In
particular, if $\k$ is perfect and $\Spec R$ is generically reduced, then it is
smooth.

\part \label{mixed char} If $\Lambda$ is torsion free and
\map{f}{\Lambda}{\Lambda_1} is an injective local map in $\hCL$
such that $\piF(\Lambda_1)$ is non-empty, then $R$ is formally smooth over $\Lambda$.
\end{Theorem}

To give the reader some feeling for this result let us quickly sketch the proof
in characteristic zero: The tangent lifting means that every horizontal tangent
vector lifts to a horizontal vector field, i.e., to a $\Lambda$-derivation. By a
result of Zariski, Lipman and Nagata this implies that $R$ is smooth.

Actually the result of these authors is more general in that they give a
structure result when only some tangent vectors lift. An analogous result would
be potentially useful also in mixed and positive characteristic and we do indeed
present such results (and we will in fact make use of the positive
characteristic version in a forthcoming article).

When trying to apply this theorem to deformations of Calabi-Yau varieties in
mixed or positive characteristic we are faced with the problem that 
the characteristic zero proof uses the
Gauss-Manin connection on de Rham cohomology to show that the de Rham cohomology
is constant
over infinitesimal thickenings. However, this constancy is in general
true only for divided power infinitesimal thickenings. To overcome this 
we shall introduce the notion of the divided power tangent lifting property that
requires the lifting property only for those $A$ that have a divided power
structure on their maximal ideals (compatible with a given such structure on
$\maxid_{\Lambda}$). In the mixed characteristic situation we again need at
least one lifting over some local $\Lambda$-algebra $\Lambda_1$ that is
{\emph{slightly ramified}}. This is defined in the next section; for the moment,
we just remark that a complete DVR of mixed characteristic $p$ is slightly ramified
over $\Z_p$ if and only if the absolute ramification index $e$ satisfies
$e<p$.

For the positive characteristic case, we define for an ideal $I$ in a ring,
$I^{(p)}$ to be the ideal generated by all $p$'th powers of $I$. We then say that a map
of local rings $(\Lambda,\maxid_\Lambda) \to (R,\maxid_R)$ with
$\maxid_\Lambda^{(p)}=0$ in characteristic $p$ is \Definition{height $1$ smooth}
if for some $n$ $R/\maxid^{(p)} \iso
\Lambda[t_1,\dots,t_n]/(t_1,\dots,t_n)^{(p)}$.
\begin{Theorem}[B]\label{Main DP Theorem}
Assume that $\piF$ (with hull $R$) fulfills the divided power tangent lifting
property.

\part If $p\Lambda = 0$ and $\maxid_\Lambda^{(p)}=0$ then $\Lambda \to R$ is
height $1$ smooth.

\part If $\Lambda$ is torsion free and \map{f}{\Lambda}{\Lambda_1} is a slightly
ramified map in $\hCL$
such that $\piF(\Lambda_1)$ is non-empty, then $R$ is smooth over $\Lambda$.
\end{Theorem}
\begin{remark} The arguments used in the proof of this force us to permit
$\Lambda$ to be a pro-artinian local ring with residue field $\k$
rather than a complete noetherian local ring. The details of these are discussed
by Gabriel, who calls them \emph{pseudo-compact} rings, in \cite{gabriel70::exp+vii+b+etude}.
\end{remark}
The theory of crystalline cohomology allows us to imitate the characteristic
zero proof of the tangent lifting property to obtain the following result where
we for simplicity have restricted $\Lambda$ to be a ring of Witt vectors.
(To repeat, in the case of characteristic zero we
merely recover a known result).
\begin{Theorem}[C]\label{Main CY Theorem}
Let $X$ be a smooth and proper purely $n$-dimensional variety over a perfect
field $\k$ with $\omega_X$ trivial. Let $R$ be a hull for the deformation functor
of $X$ either over $\k$ or, when $\ch \k > 0$ also over $\Witt=\Witt(\k)$.

\part If $\ch \k=0$ then $R$ is formally smooth.

\part If $\ch \k > 0$, $\dim_{\k} H^n_{DR}(X/\k)=\sum_{i+j=n}h^{ij}(X)$, and $R$ is
the hull for deformations over $\k$, then $R$ is height $1$ smooth.

\part If $\ch \k > 0$, $b_n(X)= H^n_{DR}(X/\k)=\sum_{i+j=n}h^{ij}(X)$ and there
exists a lifting of $X$ to a slightly ramified $\Witt$-algebra then $R$ is a
formally smooth $\Witt$-algebra.
\end{Theorem}
\begin{remark}
Schröer has proved a result that is related to the third part
\cite{schroeer01::t1}. The difference is that he requires the lifting to be done
already over $\Witt$ and that he does not specify any concrete conditions that
makes his condition on the flatness of crystalline cohomology to be fulfilled
(which we have done through the condition on Hodge and Betti numbers). Also, it
is interesting to notice that his approach is based on the $T^1$-lifting
criterion rather than tangent lifting, as here.
\end{remark}
This result may be combined with our results on the tangent lifting property as
well as some global arguments to prove the following result.
\begin{Theorem}[D]\label{CY locus}
Let $\k$ be a perfect field of characteristic $p > 0$.

\part Let \map{f}{X}{S} be an everywhere versal family of smooth and proper
purely $n$-dimensional varieties over a finite type $\k$-scheme $S$ for which
$\omega_{X_s}$ is trivial and $\dim_{\k(s)}
H^n_{DR}(X_s/\k(s))=\sum_{i+j=n}h^{ij}(X_s)$ for all $s \in S$. Then the smooth
locus of $S$ is open and closed.

\part Let \map{f}{X}{S} be an everywhere versal family of smooth and proper
purely $n$-dimensional varieties with $S$ a finite type $\Witt(\k)$-scheme for
which $\omega_{X_s}$ is trivial and $b_n(X_s/\k(s))=\sum_{i+j=n}h^{ij}(X_s)$ for
all $s \in \overline{S}$, $\overline{S} := S\Tensor\Z/p$. Then the intersection
of the smooth locus of $S$ with $\overline{S}$ equals the intersection of the
closure of $S\Tensor\Witt[1/p]$ with the smooth locus of $\overline{S}$. In
particular the smooth locus of $S$ is open and closed.
\end{Theorem}
For the reader's convenience let us point out that Theorems A and B are proved
after Proposition 3.1, Theorem C after Proposition 4.2 and Theorem D after
Proposition 4.3.
\begin{section}{Divided power envelopes and some completions}

In this section rings will not necessarily be noetherian. In particular,
``local ring'' will merely mean a ring with a unique maximal ideal.

Recall (a complete reference for this is well as other standard facts on divided
powers is \cite{berthelot74::cohom}) that a pair $(R,I)$ of a commutative ring
and an ideal in it is a \Definition{divided power pair} if there are maps
\map{\gamma_n}IR fulfilling the relations expected of maps imitating the
operations $x\mapsto x^n/n!$. For any map of pairs $(S,J) \to (R,I)$ such that $(S,J)$ is
a divided power pair we can define its \Definition{divided power hull}
$\Gamma(R,I)$ which is universal for maps of the pair
into a divided power pair whose composite with $(S,J) \to (R,I)$ is a divided
power pair. Here, $(S,J)$ is understood; for example, $(S,J)$
might be $(\Z,0)$.
Note that we always have $n!x^n=\gamma_n(x)$ so that the divided
power hull maps to the subring of $R\Tensor\Q$ generated by $R$ and the elements
$x^n/n!$ for all $x \in I$. If $R$ is the polynomial ring
$\Lambda[x_1,\dots,x_n]$ for some torsion free ring $\Lambda$ and
$I=(x_1,\dots,x_n)$ then this map is an isomorphism. By functoriality this gives
a description of the divided power hull for any pair $(R,I)$. If $J$ is a subset
contained in $I$ we define the $n$'th divided power ideal of $J$, $J^{\DP{n}}$,
as the ideal generated by the elements
$\gamma_{n_1}(x_1)\dots\gamma_{n_k}(x_k)$ for $n_1+\dots+n_k\ge n$. This is the
smallest ideal containing $J$ and stable under the $\gamma_k$. (Note that even
if $J$ is an ideal $J^{\DP1}$ is in general distinct from $J$.)

There are various topological constructions that we can now make.

\begin{definition-lemma}
\part[i] Suppose that $(R,I)$ is a DP pair where $R$ is a
topological ring, $I$ is a closed ideal and the $\gamma_k$ are 
continuous. If $J \subseteq I$ is a subset the
\Definition{$J$-DP-adic completion} of $(R,I)$ is the completion $\hat R$ in
the topology defined by closure of the divided power ideals $J^{\DP k}$. If
$\hat I$ is the image of $I$ in this completion then $(\hat R,\hat I)$ is
(continuous) divided power pair.

\part[ii] 
If $(R,\maxid_R)$ is a local ring with a divided power structure its
\Definition{pro-artinian completion} is the completion in the topology defined
by the collection of divided power ideals of finite colength contained in the
maximal ideal. (Note that to us a local ring will not necessarily be
noetherian.)

\part[iii] If $(R,\maxid_R)$ is a local ring then its
\Definition{pro-artinian DP-completion} 
is the pro-artinian completion of its
divided power hull $\Gamma(R,\maxid_R)$. It will be denoted $\hat\Gamma(R,\maxid_R)$.
This is canonically isomorphic to the inverse limit of the inverse system consisting
of finite artin $R$-algebras with a DP structure on their maximal ideals, where
the homomorphisms in the system are DP-homomorphisms. It is also canonically 
isomorphic to the pro-artinian DP-completion of the $\maxid_R$-adic completion
of $R$.

\part[iv] If $(R,\maxid_R)$ is a local ring and $I$ an ideal in $R$, then
its \Definition{$I$-DP-nilpotent pro-artinian DP-completion} 
$\hat\Gamma_I(R,\maxid_R)$ is the quotient of
$\hat\Gamma(R,\maxid_R)$ by the intersection of all the ideals $I^{\DP n}$.
This is canonically isomorphic to the inverse limit of the inverse system consisting
of finite artin $R$-algebras with a DP structure on their maximal ideals in which
the image of $I$ is DP-nilpotent, where again
the homomorphisms in the system are DP-homomorphisms. This ring is unchanged
if $(R,\maxid_R)$ is replaced by its $\maxid_R$-adic completion.
\noproof
\end{definition-lemma}

\begin{remark} \part[i] Even if we start with a noetherian ring, its
DP envelopes and their various completions defined above are not usually
noetherian. This might cast doubt on some of our arguments.  However, these
completions are pro-artinian, and the modules over them that we shall consider
will also be pro-artinian. This will allow reduction to the artinian case. For
the details of these facts, especially Nakayama's lemma, we refer to
\cite{gabriel70::exp+vii+b+etude}.

\part[ii] The pro-artinian DP-completion is not always functorial but it will be
sufficiently often, as shown by the next lemma.
\end{remark} 
\begin{lemma}\label{functorial DP complete}
A local map of local rings $f:R \to S$ for which the residue field extension is
finite extends to a continuous map between pro-artinian DP-completions. If $I,J$
are ideals of $R,S$ with $f(I)\subset J$, then the same holds for the
nilpotent DP-completions.
\begin{proof}
If $S \to A$ is a ring homomorphism where $A$ has finite length as $S$-module
whose maximal ideal has a divided power structure then it is of finite length as
an $R$-module. The second part is easy.
\end{proof}
\end{lemma}
\begin{remark}
The map $\F_p[\epsilon] \to \F_p(t)[\epsilon]$ that takes $\epsilon$ to
$\epsilon t$ and $\F_p$ to $\F_p(t)$ does not have a continuous extension to the
pro-artinian DP-completions. In fact there is a unique divided
power structure on $\F_p(t)\epsilon$ for which
$\gamma_p(\lambda\epsilon)=\lambda^p\epsilon$. The kernel of the induced map
$\Gamma(\F_p[\epsilon],\F_p[\epsilon]\epsilon) \to \F_p(t)[\epsilon]$ does not
have a kernel of finite colength, in fact $\gamma_p^n(\epsilon)$ maps to
$t^{p^n}\epsilon$ which are $\F_p$-linearly independent.
\end{remark}

\begin{definition} A local homomorphism $(\Lambda,\maxid_{\Lambda})
\to(\Lambda_1,\maxid_{\Lambda_1})$ 
is \Definition{slightly ramified} if the composite homomorphism
$\Lambda\to \hat\Gamma(\Lambda_1,\maxid_{\Lambda_1})$ is injective.
Otherwise the map is \Definition{highly ramified}.
\end{definition}

We illustrate this with some examples.

\begin{lemma} A local homomorphism $(\Lambda,\maxid_{\Lambda})
\to(\Lambda_1,\maxid_{\Lambda_1})$
of DVRs of mixed characteristic $p$ is slightly ramified if and only if
$e_1<p$, where $e_1$ is 
the absolute ramification index $e_1$ of $\Lambda_1$.
\begin{proof} This is an easy consequence of the well known fact
that $\Lambda_1$ has a DP structure on its maximal ideal  
if and only if $e_1<p$.
\end{proof}
\end{lemma}

\begin{lemma}\label{easy 0.3} Suppose that
$\Lambda\to\Lambda_1\to\Lambda_2$ are local homomorphisms of local rings
such that $\Lambda_1\to\Lambda_2$ induces a finite extension of residue
fields. If $\Lambda\to\Lambda_2$ is slightly ramified, then so is
$\Lambda\to\Lambda_1$.
\begin{proof} By (\ref{functorial DP complete}) the sequence of homomorphisms
$\Lambda\to\Lambda_1\to\Lambda_2$ maps to
a sequence 
$\Lambda\to\hat\Gamma(\Lambda_1,\maxid_{\Lambda_1})\to\hat\Gamma(\Lambda_2,\maxid_{\Lambda_2})$
where the composite $\Lambda\to\hat\Gamma(\Lambda_2,\maxid_{\Lambda_2})$ is injective.
Then $\hat\Gamma(\Lambda_2,\maxid_{\Lambda_1})$ is also injective.
\end{proof}
\end{lemma}

Now let $\Cohen$ denote a Cohen ring for the field $\k$ of characteristic $p$.

\begin{lemma} Suppose that 
the complete local domain R is finite over $\Cohen$, with ramification index $e$.
Suppose that $R$ has normalization $\sO$ and that $p = u.\pi^e$, with $\pi$
a uniformizer in $\sO$ and $u$ a unit. Say $\maxid_R.\sO=\pi^c.\sO$. Then
$\Cohen\to R$ is slightly ramified if $e<cp$.
\begin{proof} Put $S=\Cohen+\pi^c.\sO$, a subring of $\sO$ containing $R$. Then $R\to S$
is finite and it is enough to show that $\Cohen\to S$ is slightly ramified. For this, it is 
enough to show that $(S,\maxid_S)$ has a DP structure. Note that $\maxid_S=\pi^c.\sO$,
so that $(S,\maxid_S)$ has a DP structure if and only if $(\sO,\pi^c.\sO)$ does.
But this holds if and only if $v_p(\pi^c)> 1/p$, which translates as $e<cp$.
\end{proof}
\end{lemma}

\begin{corollary}\label{slice} $\pow[\Cohen]{x}/(p^a-x^b)$ is slightly ramified
over $\Cohen$ if $b<ap$.
\begin{proof} First, assume that $a,b$ are coprime. Then
$x = \pi^a$ and $p = \pi^b$. Then $e=b$ and $c=\min(a,b)$.

In general, write $a=ha_1$ and $b=hb_1$, where $h=HCF(a,b)$. Put
$R=\pow[\Cohen]{x}/(p^a-x^b)$ and $R_1=\pow[\Cohen]{x}/(p^{a_1}-x^{b_1})$. Then we have
just shown that $\Cohen\to R_1$ is slightly ramified, so that $\Cohen\to R$ is slightly
ramified by (\ref{easy 0.3}).
\end{proof}
\end{corollary}

\begin{corollary} If $\min(b,d)<ap$, then
$A:=\pow[\Cohen]{x,y}/(p^a-x^b-y^d)$ is slightly ramified over $\Cohen$.
\begin{proof} There are surjections $A\to \pow[\Cohen]{x}/(p^a-x^b)$
and $A\to \pow[\Cohen]{y}/(p^a-y^d)$. Now use (\ref{slice}) and (\ref{easy 0.3}).
\end{proof}
\end{corollary}

In particular, the $E_8$ singularities $R=\pow[\Z_2]{x,y}/(2^a-x^b-y^d)$, where
$\{a,b,d\}=\{2,3,5\}$, are slightly ramified over $\Z_2$, although in each case
every DVR that is finite over $R$ and centered at the maximal ideal is highly
ramified over $\Z_2$.

\end{section}
\begin{section}{Smooth foliations}

In this section we shall study a condition on rings that will arise as a
consequence of assuming the TLP or the DPTLP
on a deformation functor. We have in mind applications
of these idea to situations outside of those that will appear in the present
paper and hence the setting will be more general than would be needed our
current purposes.

We are going to use an argument of Zariski, Lipman and Nagata in several
slightly distinct situations and the next result is our attempt to extract the
essence of their argument in order to be able to apply it to the different
situations. Curiously enough our result has no direct relation with the
original argument and instead appears as a result on group actions (which is
essentially well known but we have not been able to find a reference which
would directly apply to our situation).
\begin{lemma}\label{map to free}
Suppose that $\cC$ is a category with fibre products, $G$ a group object in
$\cC$ and $G\times X \to X$ an action of $G$ on $X \in ob(\cC)$. Assume that we
have an equivariant map \map{s}{X}{G}, where $G$ acts on itself by left
multiplication. Define $Y$ to be the pullback of $s$ along the identity and
\map{h}{G\times Y}{X} to be the composite of the inclusion $G\times Y \to
G\times X$ with the action of $G$ on $X$.  Then $h$ is a $G$-equivariant
isomorphism, where $G$ acts on $G\times Y$ by left multiplication on the first
factor and trivially on $Y$ and $s=pr_G\circ h^{-1}$, where $pr_G:G\times Y\to
G$ is the projection.
\begin{proof}
(We shall use set-theoretic notation which can be translated into a proof using
only the given structure maps at will.) Define \map{t}{X}{Y} by
$t(x)=s(x)^{-1}x$, which maps into $Y$ as $s(t(x))=s(s(x)^{-1}x)=s(x)^{-1}\cdot
s(x)=e$. We then have the map \map{(s,t)}{X}{G\times Y} which is easily seen to
be the inverse to $h$.
\end{proof}
\end{lemma}

\begin{corollary}\label{quotient}
With the notation and assumptions of (\ref{map to free}) assume also that $\cC$
is the category of (formal) schemes.  Then \map{h^{-1}\circ pr_Y}{X}{Y} is a
quotient by the action of $G$ on $X$.
\noproof
\end{corollary}
\medskip
We shall need to generalise the notion of a smooth (height $1$) foliation from
the case of a smooth variety discussed in \cite{Ek86} to the case of singular
varieties. Just as in that case the theory is very different in the case of
characteristic zero and positive characteristic and this will soon be
apparent. Our first definition however works in any characteristic.

Let \map{f}{X}{S} be a morphism of finite type of noetherian (formal) schemes.
Assume that for every closed point of
$X$ the residue field extension of the image of the point in $S$ and that of the
point itself is separable (to us a \Definition{separable field extension} will
be algebraic).  
\begin{definition}
\part[i] A \Definition{smooth foliation} on $f$ is a subsheaf $\cF$ of
the tangent sheaf $T_{X/S}$ of the sheaf of (continuous) $S$-derivations of $X$
such that the induced map $\cF \Tensor_{\cO_X} \k(x) \to
\Hom_{\k(x)}(\maxid_x/\maxid_s+\maxid_x^2,\k(x))$ is injective for all points
$x$ of $X$, where $s:=f(x)$, and such that $\cF$ is closed under commutators.
Moreover, if $S$ is over a field of positive characteristic $p$, then we
demand that $\sF$ be closed under $p$'th powers. 

\part[ii] If $Z \hookrightarrow X$ is a closed subscheme then $\cF$
is \Definition{transverse} to $Z$ if for all points $x \in X$ the composite
$\cF \Tensor_{\cO_X} \k(x) \to \Hom_{\k(x)}(\maxid_x/\maxid_s+\maxid_x^2,\k(x)) \to
\Hom_{\k(x)}(I/\maxid_xI,\k(x))$ is injective, where $I$ is the ideal of $Z$.

\part [iii] An $S$-map \map{f}{X}{Y} is a \Definition{quotient} of $\sF$ if $\sF$
is the relative tangent sheaf of it and $\sO_Y$ consists of the $\sF$-constants
(the sections annihilated by the local sections of $\sF$) of $f_*\sO_X$.
\begin{remark}
\part When $f$ is smooth, the injectivity condition simply says that $\cF$ is a
subbundle of $T_{X/S}$ so this corresponds to the already established notion.

\part Note that a smooth foliation is transverse to every closed point.

\part It is not clear that our definition is the right one
in the mixed characteristic case. The reason is that we do not
know what to do with the condition of being closed under $p$'th powers. On the
one hand it makes no sense in characteristic zero, on the other hand in the
mixed characteristic situation, the reduction modulo $p$ of a smooth foliation
should be a smooth foliation and hence some condition is needed. One could
require that for $D$ in the foliation $D^p$ should be of the form $E+pF$ for $p$
a prime, $F$ a differential operator and $E$ in the foliation but that seems
rather artificial.

A similar problem occurs when later in this section we shall consider divided
power foliations. Again, if the divided power structure is trivial, i.e., when
the divided power ideal is the zero ideal one needs some $p$-integrability
condition. On the other hand, even in positive characteristic $p$ the
$p$-integrability for divided power derivations does not make sense as the
$p$'th power of a divided power derivation is not a divided power
derivation. However, our definition turns out to be adequate for the purposes
of this paper.
\end{remark}
\end{definition}
We shall see that, just as in the smooth case, a smooth foliation gives rise
to a flat infinitesimal equivalence relation on $X$. We start with a
preliminary result on modules.

Suppose that $(A,\maxid,\k)$ be a local ring and that $M,N$ are $A$-modules
with $N$ a finitely generated. Assume given
an $A$-homomorphism $\phi:N \to M^* := \Hom_A(M,A)$ and denote by $\phi^*$
the composite homomorphism $M\to M^{**}\to N^*$.

\begin{lemma}\label{direct factor}
Assume 
that either the map $\bar\phi=\phi\Tensor 1:N/\maxid N \to \Hom_\k(M/\maxid M,\k)$ 
induced by $\phi$ is injective or $N = M^*$ and $\bar\phi$ is surjective. 
Then the following are true.

\part[i] $\phi$ is injective and $N$ is free.

\part[ii] $M$ is a direct
sum $M_1 \bigoplus M_2$ such that $M_1^*$ is identified with $\phi(N)$,
$N$ and $M_2$ are mutual annihilators and $N/\maxid N$ and 
$M_2/\maxid M_2$ are mutual annihilators.

\part[iii] Given elements
$m_1,\dots,m_n$ of $M_1$, they form a basis for $M_1$ if and only if
the elements 
$\phi^*(m_1),\ldots,\phi^*(m_n)$ form a basis of $N^*$, or if and only if
$\phi^*(m_1),\ldots,\phi^*(m_n)$ map to elements forming
a basis of $\Hom_\k(N/\maxid N,\k)$ under the natural homomorphisms
$N^*\to \Hom_\k(N/\maxid N,\k)$.
\begin{proof}
Assume first that $\bar\phi$ is injective.  Let 
$n := \dim_\k N/\maxid N$ and pick elements $m_1,\dots,m_n$ of $M$ 
such that the composite of 
$\bar\phi$ and the map $\Hom_\k(M/\maxid M,\k) \to \k^n$ 
given by evaluation at
the residues $\overline{m_i} \in M/\maxid M$ is an isomorphism. This gives us a map
$A^n \to M$ given by $(r_i) \mapsto \sum_ir_im_i$.  Evaluation at the $m_i$
gives its transpose $M^* \to A^n$ and the composite $N \to M^* \to A^n$ induces
an isomorphism upon reduction modulo $\maxid $. By Nakayama's lemma it is then
surjective and by the projectivity of $A^n$ and Nakayama's lemma again it is an
isomorphism. This shows that $N$ is a direct factor of $M^*$ as well as
isomorphic to $A^n$. By construction the basis $(n_i)$ of $N$ thus constructed
is dual to $\{m_i\}$ so that the composite $A^n \mapright{(m_i)} M
\mapright{(n_i)} A^n$ is the identity and hence the $m_i$ generate a direct
summand of $M$. This shows that $M$ is isomorphic to $M_1\Dsum M_2$ in such a
way that $N \to M^*$ is the transpose of the projection $M_1\Dsum M_2 \to
M_1$. Using this description the second statement is clear and the last is
obvious.

If instead $N=M^*$ and $\bar\phi$ is surjective then we can
choose $m_1,\dots,m_n \in M$ which form a basis for $M/\maxid M$ and then by
assumption find $f_1,\dots,f_n \in M^*$ such that $f_im_j \equiv \delta_{ij}
\bmod \maxid $. This gives a map $M \to A^n$ which is surjective by Nakayama's lemma
and then an isomorphism as it is split. In particular we have the first
situation and the results already obtained apply.
\end{proof}
\end{lemma}
An immediate consequence of this lemma is that
the tangent sheaf is a smooth foliation under conditions weaker than those
of the definition.
\begin{proposition}\label{tangent lifting}
Let \map{f}{X}{S} be a morphism of finite type of noetherian (formal) schemes.
If for every point $x$ of $X$ the map $T_{X/S,x} \to
\Hom_{\k(x)}(\maxid_x/\maxid_{f(x)}+\maxid_x^2,\k(x))$ is
surjective, then $T_{X/S}$ is a smooth foliation.
\begin{proof}
By (\ref{direct factor}) the condition implies that the natural evaluation homomorphism
$T_{X/S,x}\Tensor\k(x)
\to \Hom_{\k(x)}(\maxid_x/\maxid_{f(x)}+\maxid_x^2,\k(x))$ is an isomorphism.
\end{proof}
\end{proposition}

\begin{definition} The map $f$ is \Definition{tangent smooth} 
when $T_{X/S}$ is a smooth foliation.
\end{definition}

We now show that smooth foliations
on a singular variety behave much like smooth foliations on a smooth variety.
\begin{theorem}
Let \map{f}{X}{S} be a morphism of finite type of (formal) schemes with $S$
noetherian and $\sF$ be a smooth foliation for $f$ transverse
to a closed subscheme $Z$ of $X$.

\part[direct] \label{foliation splitting}
Locally $\Omega_{X/S}^1$ splits as a sum $\sG'\Dsum\sG''$, where $\sG'$ is
free with a basis of the form $dx_i$, $\sF=(\sG')^*$ and $\sF$ is identified with the
annihilator of $\sG''$. In particular, $\sF$ is locally free and locally
$T_{X/S}$ splits as $\sF\Dsum (\sG'')^*$. Furthermore, 
at a point of
$Z$ the $x_i$ can be chosen to lie in $I_Z$.

\part[char 0] If $S$ is over $\Sp\Q$ let $\hat X$ be the formal completion of $X$ along
$Z$ and $\hat S$ that of $S$ along the image of $Z$.
Then the restriction of $\sF$ has a quotient $\hat X \to Y$ which is formally
smooth and the composite $Z \to \hat X \to Y$ is an embedding. If $\sF=T_{X/S}$
(which is possible precisely when $f$ is tangent smooth) then the structure map
$Y \to \hat S$ is an immersion.

\part[char p] If $S$ is over $\Sp\F_p$ then $\sF$ has a quotient $X \to Y$ which is flat
and the composite $Z \to X \to Y$ is an embedding. Locally $X \to Y$ has the
form $\Sp \sO_Y[t_1,\dots,t_n]/(t_1^p-f_1,\dots,t_n^p-f_n)$ with the $t_i$
vanishing on $Z$.
\begin{proof}
For the first part, Lemma \ref{direct factor} shows that
$\Omega^1_{X/S}=M_1\Dsum\Ann(\sF)$ with $M_1=\sF^*$. 
Moreover, for all $x\in Z$ there exist
$x_1,\ldots,x_n\in I_Z$ and $D_1,\ldots,D_n\in\sF$ such that $D_1,\ldots,D_n$
give a basis of $\sF\Tensor\k(x)$ and the matrix
$(D_i(x_j)\pmod{\maxid_x})$ is an invertible
$n\times n$ matrix over $\k(x)$. This completes the proof of
(\ref{foliation splitting}).

Suppose that we are over $\Sp\Q$. We use
a slight modification of the proof of
a result of Zariski, Lipman and Nagata (\cite[p.~230]{matsumura86::commut}). Since
quotients are unique, when they exist, it is enough to
work locally on $Z$.  By assumption and \DHrefpart{direct} we can, on an
affine neighbourhood $\Sp A$ of a point $x \in Z$, find a basis $D_1,\dots,D_n$
of $\sF$ and functions $x_1,\dots,x_n$ vanishing on $Z$ such that $D_ix_j =
\delta_{ij}$. This implies that $[D_i,D_j]x_k=0$, so that
$[D_i,D_j]\in\Ann(\sG')$. Since the vector fields $[D_i,D_j]$ lie in 
$\sF$, by assumption, and $\sF=\Ann(\sG'')$, it follows that
$[D_i,D_j]=0$ for
all $i$ and $j$.  

We write $\hX=\Spf R$
and $\hS=\Spf\Lambda$. Let $\mathbf{\hat G}$ denote the formal group
$\mathbf{\hG}_a^n=\Spf \pow[\Lambda]{t_1,\ldots,t_n}$, with formal co-multiplication
$t_i\mapsto t_i\hat\otimes 1 + 1\hat\otimes t_i$. Then the
(continuous extensions of the) derivations $D_i$ define an action of the formal
group $\mathbf{\hG}$ on $\Spf R$ by the
formal co-action
\begin{displaymath}
r \mapsto \sum_{\alpha \in \N^n}\frac{t^\alpha}{\alpha!}D^\alpha(r),
\end{displaymath}
That this is a
co-action follows from the fact that $[D_i,D_j]=0$. We also define a map
\map{s}{\Spf R}{\mathbf{\hG}_a^n} given by $t_i \mapsto x_i$. The condition
that $D_ix_j=\delta_{ij}$ is then precisely that $s$ is equivariant and we may
now apply (\ref{map to free}). 

Finally, assume that $\sF=T_{X/S}$.
Then the surjectivity of $\Lambda \to C$, 
where $C\subset R$ is the ring of $\sF$-invariants,
follows because modulo $(x_1,\dots,x_n)$ it is surjective.

The characteristic $p$ case is done similarly. We work locally on $X$ and
get $D_i$ and $x_j$ in the
same way, as well as the relation $[D_i,D_j]=0$. Moreover,
$D_i^p$ is a linear combination of the $D_j$ which gives the relation $D_i^p=0$.
The formula
\begin{displaymath}
r \mapsto \sum_{0\le \alpha_i <p}\frac{t^\alpha}{\alpha!}D^\alpha(r),
\end{displaymath}
where the $t_i$ are parameters on $\alpha_p^n=\Sp
\k[t_1,\dots,t_n]/(t_1^p,\dots,t_n^p)$, defines an action of
$\alpha_p^n$ on $X$. This time however the $x_i$ do not
necessarily give a map to $\alpha_p^n$ as we may not have $x_i^p=0$. Instead we
consider the $\alpha_p^n$-torsor $\Sp \k[r_1,\dots,r_n] \to \Sp \k[s_1,\dots,s_n]$
given by $s_i \mapsto r_i^p$ and the Lie algebra of $\alpha_p^n$ acting as
$\partial/\partial r_i$. Then $X \to \Sp \sO_S[r_1,\dots,r_n]$ given by $r_i \mapsto
x_i$ is an equivariant map and as $\alpha_p^n$ acts freely on $\Sp
\sO_S[r_1,\dots,r_n]$ it does so on $X$. As $\alpha_p^n$ is a finite group scheme
there is no problem with taking the quotient and we have a quotient map $X \to
Y$ and an induced map $Y \to \Sp \sO_S[s_1,\dots,s_n]$. Being a map between group
torsors this mean that $X$ is the fibre product of $Y \to \Sp \sO_S[s_1,\dots,s_n]$
and $\Sp \sO_S[r_1,\dots,r_n] \to \Sp \sO_S[s_1,\dots,s_n]$ which gives the desired
result. Again the claim made about $Z$ is clear.
\end{proof}
\end{theorem}
We can apply the theorem to the case where $T_{X/S}$ is a smooth foliation. In
that case we also get a result in the mixed characteristic case. 


\begin{corollary}\label{T-smooth}
Let \map{f}{X}{S} be a dominant morphism of finite type of (formal) schemes with
$S$ noetherian.

\part[0] If $f$ is smooth then $f$ is tangent smooth.

\part[i] $\Omega^1_{X/S}$ is locally free if and only if $f$ is tangent smooth.

\part[ii] If $S$ is over $\Sp\Q$ and $f$ is tangent smooth then $f$ is smooth.

\part[iii] If $S$ is over $\F_p$ then $f$ is tangent smooth if and only if
for every point $x \in X$ and $s:=f(x)$ there is an isomorphism
$\widehat{\sO}_{X,x} \iso \pow[\widehat{\sO}_{S,s}]{t_1,\dots,t_n}/(f_i)$, where
$f_i(t_1,\dots,t_n)
= g_i(t_1^p,\dots,t_n^p)$ for some $g_i \in \pow[\widehat{\sO}_{S,s}]{t_1,\dots,t_n}$.

\part[iv] If $\sO_S$ is torsion-free and $f$ is tangent smooth then $f$ is smooth.
\begin{proof}
\DHrefpart{ii} and \DHrefpart{iii} follow directly from the theorem and
(\ref{tangent lifting}) together (in the characteristic zero case) with the fact that
a dominant immersion is \'etale.

For the last part, we may localize $X$ and $S$ and then complete, to
get $X=\Spf R$ and $S=\Spf\Lambda$. Then choose
$t_1,\dots,t_n$ in the maximal ideal $\maxid_R$ of $R$ that map to a basis
of the cotangent space $t_{R/\Lambda}^*=\maxid_R/(\maxid_R^2+\maxid_\Lambda)$
and define a map $\sO:=\pow[\Lambda]{T_1,\dots,T_n} \stackrel{\pi}{\to}R$ with
$\pi(T_i)= t_i$. Then $\pi$ is surjective.

Since $\Lambda$ is torsion-free,
there is an embedding $\Lambda\to K$, where $K$
is a finite direct product of characteristic zero fields.
There is a commutative diagram
$$\begin{array}{cccc}
\sO & \stackrel{\pi}{\longrightarrow} & R\\
\downarrow & & \downarrow\\
\pow[K]{T_1,\dots,T_n}&\stackrel{\pi_1}{\longrightarrow} & R{\widehat\otimes}_\Lambda K,
\end{array}$$
where $R{\widehat\otimes}_\Lambda K$ is the tensor product completed
in the $(t_1,\dots,t_n)$-adic topology.
Since $\pi_1$ is an isomorphism, by \DHrefpart{ii}, and $\sO\to\pow[K]{T_1,\dots,T_n}$
is injective, it follows that $\pi$ is injective.
\end{proof}
\end{corollary}
We have put the positive characteristic case alongside the characteristic zero
and mixed characteristic case but in actuality it is rather different and we
finish this section with a discussion special to that case.

We shall say that a local ring $(R,\maxid)$ with $pR=0$, $p$ a prime, is of
\Definition{height $1$} if $f^p=0$ for all $f \in \maxid$. We shall also use the
notation $\maxid^{(p)}$ for the ideal generated by the $p$'th powers of elements
of $\maxid$. It is then clear that $R/\maxid^{(p)}$ is the largest height $1$
quotient of $R$. We shall also say that $R$ is \Definition{formally height
$1$-smooth} if it has the infinitesimal lifting property with respect to maps
into height $1$ local rings.
\begin{proposition}\label{height 1 characterisation}
Let $p$ be a prime and $(R,\maxid)$ a local noetherian ring with $pR=0$ and
perfect residue field $\k$. The following conditions are equivalent:
\begin{enumerate}
\renewcommand\labelenumi{(\roman{enumi})}
\item $R$ is formally height $1$-smooth.

\item $R/\maxid^{(p)}$ is of the form
$\pow[\k]{t_1,\dots,t_n}/(t_1,\dots,t_n)^{(p)}$ for some $n$.
\end{enumerate}
\begin{proof}
This is completely analogous to the similar relation between the infinitesimal
lifting property and formal smoothness:

We may assume that $R$ is complete and then, as $\k$ is perfect, $R$ contains $\k$
and is thus the quotient of $\pow[\k]{t_1,\dots,t_n}$ with $n$ minimal. We then
use the height $1$ lifting property for the map
$\pow[\k]{t_1,\dots,t_n}/(t_1,\dots,t_n)^{(p)} \to
\pow[\k]{t_1,\dots,t_n}/(t_1,\dots,t_n)^2$ to get a map $R \to
\pow[\k]{t_1,\dots,t_n}/(t_1,\dots,t_n)^{(p)}$. This induces a map
$R/\maxid^{(p)} \to \pow[\k]{t_1,\dots,t_n}/(t_1,\dots,t_n)^{(p)}$ whose
composite with $\pow[\k]{t_1,\dots,t_n}/(t_1,\dots,t_n)^{(p)} \to
R/\maxid^{(p)}$ is an automorphism of $\pow[\k]{t_1,\dots,t_n}/(t_1,\dots,t_n)^{(p)}$.
\end{proof}
\end{proposition}

As mentioned in the introduction, when applying our results to deformation of
Calabi-Yau varieties we shall want to restrict ourselves to deformations over
elements of $C_{\Lambda}$ which have a divided power structure on their maximal
ideals. As should perhaps be no surprise the results of this section become more
uniform (and hence more like the results in characteristic $0$) if we modify our
arguments by using divided powers. 

We need to recall some facts about continuous derivations and differentials
for continuous maps of
topological divided power pairs $(\Lambda,J) \to (R,I)$.
A (continuous) divided power derivation consists of a (topological)
module $M$ over $R$ and a (continuous)
$\Lambda$-linear map \map{D}{R}{M} such that $D(rs)=rD(s)+sD(r)$ and
$D\gamma_n(i)=\gamma_{n-1}(i)D(i)$ for all $r,s \in R$ and $i \in I$. There then
is a universal (continuous) divided power derivation
\map{d}{R}{\Omega^{1,DP}_{R/\Lambda}}. It has the property that if $\maxid_R
\supseteq I$ is a maximal ideal that induces a maximal ideal $\maxid_\Lambda$ of
$\Lambda$ such that the residue field extension is separable then
$\Omega^{1,DP}_{R/\Lambda}\Tensor R/\maxid_R$ is (canonically) isomorphic to
$\maxid_R/\maxid_R^2+I^{\DP2}+\maxid_\Lambda$. 

\begin{definition} A map $(\Lambda,J) \to (R,I)$ as above
is \Definition{of finite type} if it
arises from taking DP envelopes, or the pro-artinian completions of these,
of an essentially finite type homomorphism.
\end{definition}

\begin{lemma} If $(\Lambda,J) \to (R,I)$ is a continuous map
of topological divided power pairs, then $\Omega^{1,DP}_{R/\Lambda}$
is a finitely generated $R$-module.
\noproof
\end{lemma}

For a (continuous) map of (topological) divided power pairs $(\Lambda,J) \to (R,I)$ we let
$\Der^{DP}_{R/\Lambda}$, which we shall also denote $T_{R/\Lambda}^{DP}$,
consist of the (continuous) divided power derivations of $R$. The pair
$(R[\epsilon],I[\epsilon])$ has a unique (topological) divided power structure extending that
of $(R,I)$ and having $\gamma_n(\epsilon)=0$ if $n > 1$. It is then easy to see
that a (continuous) divided power map $R \to R[\epsilon]$ which is the identity modulo
$\epsilon$ is the same thing as a (continuous) divided power derivation of $R$.

There is now an obvious extension of the notions of a smooth foliation 
and a tangent smooth map to the divided power case.
We express this in a global context.

\begin{definition}
Suppose that $(X,Y) \to (S,T)$ is a divided power map of divided power pairs of (formal)
schemes and that the residue fields of all closed points are perfect. 

\part A \Definition{smooth DP foliation} (transverse to $Y$) on the map is an
$\sO_X$-submodule $\sF$ of $T_{X/S}^{DP}$ which is closed under commutators
such that for all closed points $x\in X \setminus Y$ the natural map
$\sF/\maxid_x \sF \to
\Hom_{\k(x)}(\maxid_x/\maxid_x^2+I^{\DP2}+\maxid_{f(x)},\k(x))$ is injective  
and for each closed point 
$y \in Y$ the induced map $\sF\Tensor \k(y) \to \Hom_{\k(y)}(I_Y\Tensor \k(y),\k(y))$ 
is injective.

\part The map is \Definition{DP tangent smooth} if $T_{X/S}^{DP}$ is a smooth DP foliation.
\end{definition}

%

There are immediate analogues of Proposition \ref{tangent lifting} and Theorem
\ref{foliation splitting}. We give these together in the next lemma.  The proofs
are identical, so omitted.

\begin{lemma} 
Suppose that \map{f}{(X,Z)}{(S,T)} be a finite type divided power morphism
of (formal) divided power schemes.

\part[global]\label{global}
If for any point $x$ of $X$, the map $T_{X/S,x}^{DP} \to
\Hom_{\k(x)}(\maxid_x/(\maxid_{f(x)}+\maxid_x^2+I_Z^{\DP 2}),\k(x))$ is surjective, 
then $T_{X/S}$ is a smooth DP foliation.

\part[DP splitting] Suppose that $\sF$ is a smooth DP foliation for $f$ transverse
to $Z$. Then, locally on $X$, we have $\Omega^{1,DP}_{X/S}=\sG'\Dsum\sG''$,
where $\sG'$ is free with a basis of the form $\{dx_i\}$, $\sF=(\sG')^*$
and $\sF$ is identified with the annihilator of $\sG''$. In particular,
$\sF$ is locally free and locally $T^{DP}_{X/S}$ splits as $\sF\Dsum(\sG'')^*$. 
Furthermore, at a point of $Z$ the $x_i$ can be chosen to lie in $I_Z$.
\noproof
\end{lemma}

Now fix a complete DP local ring $\Lambda$
and consider the category $\cC$ of local pro-DP-artin $\Lambda$-algebras $R$
with a section $R\to\Lambda$ and the further property that $\ker(R\to\Lambda)$
is (topologically) DP-nilpotent. This category has finite sums, denoted by
$\hot$, whose construction is left to the reader.

\begin{lemma}\label{compatibility} If $A\in\cC$, the 
DP envelope $\Gamma$ of the pair $(A[t_1,\ldots,t_n],(t_1,\ldots,t_n))$ has a DP
structure on the maximal ideal $\maxid$ generated, as a DP ideal, by
$(\maxid_A,t_1,\ldots,t_n)$ and the ring $A\la\la t_1,\ldots,t_n\ra\ra$
constructed as the pro-DP-artin completion of $\Gamma$ at $\maxid$ lies in
$\cC$.
\begin{proof} Denote by $M$ the free $A$-module with basis $\{t_1,\ldots,t_n\}$.
Then $\Gamma=\bigoplus \Gamma^n(M)$ and 
$A\la\la t_1,\ldots,t_n\ra\ra =\prod \Gamma^n(M)$.
Now the result, in particular the fact that $A\la\la t_1,\ldots,t_n\ra\ra$ 
is $(t_1,\ldots,t_n)$-DP-adically separated, is clear.
\end{proof}
\end{lemma}

\begin{theorem}
\part[i] Suppose that $\sF$ is a smooth DP foliation for the object
$\Lambda\to R$ of $\cC$, transverse to the given section $\Spf\Lambda$ of $\Spf R$.
Put $\Lambda '=R^{\sF}$, the ring of $\sF$-invariants.
Then $\Lambda '$ is in $\cC$
and $R$ is isomorphic to 
$\Lambda '\la\la t_1,\ldots,t_n\ra\ra$.
\part[ht1]\label{DP to ht 1}
Assume that $\k$ is a perfect field of positive characteristic $p$ and let
$(\Lambda,\maxid_\Lambda)$ be a local $\k$-algebra whose residue field is a
finite extension of $\k$. Assume that its divided power
hull is DP tangent smooth over $\k$. Then it is height $1$ smooth.

\part\label{DP to P} Assume that $\Lambda_0 \stackrel{f}{\to} R\to \Lambda_1$ 
are local maps of noetherian local rings,
that the induced extensions of residue fields are finite, 
that $\Lambda_0$ is torsion-free and that $\Lambda_0 \to\Lambda_1$
is slightly ramified.
Suppose that
$\hat f:\widehat\Lambda_0 \to\widehat R$ is the induced homomorphism between
the pro-DP-artin completions of the DP envelopes of the first two rings.
If $\hat f$ is DP tangent smooth, then $f$ is 
formally smooth.

\begin{proof}
The proof of the first statement is essentially identical to what has already
been proven in characteristic zero; we get commuting generators $D_1,\ldots,D_n$
of $\sF$ and topological DP generators $x_1,\ldots,x_n$ of $\ker(R\to\Lambda)$
such that $D_i(x_j)=\delta_{ij}$. The difference is that the formal group
${\hG}_a^n$ must be replaced by the commutative formal group ${\hG}=\Spf{\hGa}$,
where $\hGa =\Lambda\la\la t_1,\ldots,t_n\ra\ra$ and the formal
co-multiplication is $t_i\mapsto t_i\hot 1 + 1\hot t_i$.  There is a formal
co-action $\hR\to \hR\ {\hot}\ {\hGa}$ given by
$$r\mapsto \sum D^\alpha(r)\ {\hat\otimes}\ \gamma_\alpha(t),$$ 
where $\alpha=(\alpha_1,\ldots,\alpha_n)$ and 
$\gamma_\alpha(t)=\gamma_{\alpha_1}(t_1)\cdots\gamma_{\alpha_n}(t_n)$.
(The $(t_1,\ldots,t_n)$-DP-adic separation is needed to ensure 
the convergence of this formula.)
As in characteristic zero, there is, because $\hR$ is
$(x)$-DP-adically separated, a $\hG$-equivariant map
\map{s}{\Spf \hR}{\hG} given by $t_i\mapsto x_i$, 
and we conclude by (\ref{map to free}) and (\ref{quotient}).


As for \DHrefpart{ht1} as the residue field of $\Lambda$ is also perfect we may
assume that it equals $\k$. we may choose elements basis $t_1,\dots,t_n \in
\maxid_{\Lambda}$ such that they form a basis for
$\maxid_{\Lambda}/\maxid_{\Lambda}^2$. We then have a surjective ring
homomorphism $\pow[\k]{t_1,\dots,t_n} \to \Lambda$. Let $R$ be the divided power
hull of $\Lambda$. Using the notations of the first part we then have $R \riso
\Lambda '\la\la t_1,\dots,t_n\ra\ra$ and we have a $\k$-homomorphism $\Lambda '
\to k$. Composing we get maps $\pow[\k]{t_1,\dots,t_n}/\maxid^{(p)} \to
\Lambda/\maxid_\Lambda^{(p)} \to \k\la\la t_1,\dots,t_n\ra\ra/\maxid^{[p]}$ and
as the composite is an isomorphism, the first map is injective and as it is also
surjective it is an isomorphism. We then conclude by Proposition \ref{height 1
characterisation}.

For the last part, we can assume that the rings are complete.
We let $\widehat\Lambda_i$ denote the pro-DP-artin completion of the DP hull
of $(\Lambda_i,\maxid_{\Lambda_i})$. Moreover, after replacing $\Lambda_1$ by its image
in $\widehat\Lambda_1$, we can assume that $\Lambda_1$ injects into $\widehat\Lambda_1$.
Of course, $\Lambda_0$ injects into $\Lambda_1$.

Choose $x_1,\ldots,x_n\in\maxid_R$ that map to a basis of the cotangent space
$\maxid_R/(\maxid_R^2+\maxid_{\Lambda_0})$. Then there is a surjection
$\pi:\sO_0:=\pow[\Lambda_0]{t_1,\ldots,t_n}\to R$ with $t_i\mapsto x_i$.
Say $J=\ker\pi$; then $J\subset \maxid_{\sO_0}^2+\maxid_{\Lambda_0}\sO_0$.

Put $R_1=R\ {\widehat\otimes}_{\Lambda_0}\ \Lambda_1$ and 
$\sO_1=\sO_0\ {\widehat\otimes}_{\Lambda_0}\ \Lambda_1$.
Write $x_i$ for $x_i \ {\widehat\otimes}\ 1$ and $t_i$ for $t_1 \ {\widehat\otimes}\ 1$.
Let $\widehat\sO_1$, resp., $\widehat R_1$, be the $(t)$-DP-nilpotent,
resp. $(x)$-DP-nilpotent, pro-artin DP-completion of $\sO_1$, resp. $R_1$.
Then the sequence $\Lambda_1\to \sO_1\to R_1\to\Lambda_1$ maps to a sequence
$\widehat\Lambda_1\to \widehat\sO_1\to \widehat R_1\to\widehat\Lambda_1$.
Then \DHrefpart{i} (more precisely, its proof) shows that 
${\hat\pi}_1:\widehat\sO_1\to \widehat R_1$
is an isomorphism.

Consider the square
$$\begin{array}{clcr}
\sO_1 & \stackrel{\pi_1}{\longrightarrow} & R_1\\
\downarrow & & \downarrow\\
\widehat \sO_1 & \stackrel{{\hat\pi}_1}{\longrightarrow} & \widehat R_1.
\end{array}
$$
Now $\sO_1$ is, as a $\Lambda_1$-module, the direct product $\prod S^m(M)$,
where $M$ is the free $\Lambda_1$-module with basis $\{t_1,\ldots,t_n\}$,
while $\widehat\sO_1$ is, as a $\widehat\Lambda_1$-module, the direct product
$\prod \Gamma^m(M\otimes_{\Lambda_1}\widehat\Lambda_1)$. Since $\Lambda_0$
is torsion-free, the natural map
$\sO_1\to\widehat\sO_1$ given by $t^m\mapsto m!\gamma_m(t)$ is injective. 
Since ${\hat\pi}_1$ is an isomorphism,
it follows that $\pi_1$ is injective, so an isomorphism.

Next, consider the commutative diagram with exact rows
$$\begin{array}{cccccccccc}
0 & \to & J & \to & \sO_0 & \stackrel{\pi}{\longrightarrow} & R &\to & 0\\
& & \downarrow & & \downarrow & & \downarrow \\
0 & \to & J_1 & \to & \sO_1 & \stackrel{\pi_1}{\longrightarrow} & R_1 & \to & 0.
\end{array}
$$
Since $\pi_1$ is injective and $\Lambda_0\to \Lambda_1$ is injective,
$\sO_0\to\sO_1$ is injective, and so, by the snake lemma, $J=0$.
\end{proof}
\end{theorem}

\end{section}
\begin{section}{The tangent lifting property}
In preparation for
the study of deformations of Calabi-Yau manifolds we shall say that 
$\piF$ fulfills the
\Definition{divided power tangent lifting property} (abbreviated to DPTLP) if
\begin{displaymath}
\piF(A[\epsilon]) \to \piF(A) \times \piF(\k[\epsilon])
\end{displaymath}
is surjective for all $A \in C_{\Lambda}$ for which the pair $(A,\maxid_A)$
admits a divided power structure. Note that under condition $H_2$ on $\piF$ this
is equivalent to $\piF(A[\epsilon]) \to \piF(A\times_{\k}\k[\epsilon])$ being surjective.
\begin{proposition}\label{TLP}
Let $\Lambda$ be complete local ring with residue field $\k$, a perfect field,
and assume that \map \piF{C_{\Lambda}}{\underbar{Sets}} is a covariant functor from
the category of local artinian $\Lambda$-algebras with residue field $\k$ to the
category sets such that $\piF(\k)$ is a point (we call this condition
$H_0$) and $\piF$ fulfills the conditions ($H_1$-$H_3$)
of Schlessinger (cf.\ \cite[Thm.\ 2.11]{schlessinger68::funct+artin}). Let $R$ be a hull for
$\piF$.

\part[i] Let \map{\pi}{\piH}{\piG} be a formally smooth map between functors
$C_{\Lambda} \to \underbar{Sets}$ that
both fulfill conditions $H_0$ and $H_2$. Then $\piH$ fulfills the (divided
power) tangent lifting property precisely when $\piG$ does. In particular $\piF$
fulfills the (divided power) tangent lifting property precisely when $h_R$ does
(where $h_R(A)=\Hom_\Lambda(R,A)$).

\part[ia] $h_R$ fulfills the tangent lifting property precisely
when $R$ is tangent smooth.

\part $h_R$ fulfills the divided power tangent lifting property
precisely when $\hat\Gamma(R,\maxid_R)$ is divided power tangent smooth.
\begin{remark}
It will be clear from the proof that other conditions than $H_2$ would also give
\DHrefpart{i}, for instance that $\piH(\k[\epsilon]) \to \piG(\k[\epsilon])$ is a bijection.
\end{remark}
\begin{proof}
Assume that $\piG$ has the tangent lifting property. Suppose $A\in C_\Lambda$,
and $\phi\in \piH(A\times_{\k}\k[\epsilon])$. We must, by the comment made at the
beginning of this section, find $\Phi\in \piH(A[\epsilon])$ with $\Phi\mapsto
\phi$. Let $\xi \in \piG(A\times_{k}\k[\epsilon])$ be the image under $\pi$ of
$\phi$. By the tangent lifting property there is a $\Xi \in \piG(A[\epsilon])$
mapping to $\xi$. By the formal smoothness of $\pi$ and the fact that
$A[\epsilon] \to A\times_{\k}\k[\epsilon]$ is surjective there is a $\Phi \in
\piH(A[\epsilon])$ mapping to $\Xi$ and $\phi$. The converse is similar (but
easier). The last statement of \DHrefpart{i} follows from the fact that as $R$
is a hull there is a formally smooth map $h_R \to \piF$.

Turning to \DHrefpart{ia} assume that $h_R$ has the tangent lifting
property. Write $R = \ili R_\alpha$ as a limit of local artinian rings with
surjective transition maps inducing an isomorphisms on cotangent spaces. Given $v
\in t_R$ and an index $\alpha$ there is a $\Lambda$-derivation
\map{D}{R}{R_\alpha} inducing the given $v$ on cotangent spaces. The set of such
derivations is a coset for the kernel of the map
$\Hom_R(\Omega^1_{R/\Lambda},R_\alpha) \to
\Hom_{\k}(\maxid_R/\maxid_R^2+\maxid_\Lambda,\k)$. That kernel is of finite length
as $R$-module so by linear compactness there is an element $D$ of
$\ili\Hom_R(\Omega^1_{R/\Lambda},R_\alpha)$ mapping to $v$. The $D$ corresponds
to an inverse system of derivations \map{D_\alpha}{R}{R_\alpha} and passing to
the limit gives a derivation $R \to R$. The converse is clear.

In the divided power case we let $(\hat R,\hat\maxid) \in \hCL$ be the completed
divided power hull of $(R,\maxid)$. By definition we have $\hat R = \ili
R_\alpha$, where $R_\alpha$ form an inverse system of divided power local
artinian algebras (and divided power maps between) inducing an isomorphisms on
(divided power) cotangent spaces. Fix a $v \in t^{DP}_{\hat R}=t_R$. For any
given $\alpha$ we may consider the composite $R \to \hat R \to R_\alpha$ and by
the divided power tangent lifting property there is a ring homomorphism $R \to
R_\alpha[\epsilon]$ inducing $v$ on cotangent spaces. We may give
$R_\alpha[\epsilon]$ the divided power structure which is the given one on the
maximal ideal of $R_\alpha$ and for which all higher divided powers of
$\epsilon$ are zero. By the universal property of $\hat R$ this extends to a
divided power map
$\hat R \to R_\alpha[\epsilon]$ reducing modulo $\epsilon$ to the given
one. This then corresponds to a divided power derivation $\hat R \to R_\alpha$
inducing $v$ on cotangent spaces. The rest of the argument is then the same as
for \DHrefpart{ia}.
\end{proof}
\end{proposition}
Combining this result with those of the previous section we get Theorem \ref{Main
Theorem}.
\begin{proofof}{Theorem \ref{Main Theorem}}
The result follows from (\ref{TLP}) and Corollary \ref{T-smooth} together with the
observations that in the positive characteristic case, if some $f_i \ne 0$ then
it is a non-zero nilpotent element of $\pow[\k]{t_1,\dots,t_n}/(f_i^p)$ whose
support is all of $\Sp \pow[\k]{t_1,\dots,t_n}/(f_i^p)$ and that in the mixed
characteristic case the existence of $\Lambda_1$ forces $\Spf R \to \Spf \Witt$ to
be dominant.
\end{proofof}
\begin{proofof}{Theorem \ref{Main DP Theorem}}
This is a direct corollary of Proposition \ref{TLP}.
\end{proofof}
\end{section}
\begin{section}{Deformation of Calabi-Yau varieties}

The results of the previous section can now be applied to the deformation of
varieties and in particular to Calabi-Yau varieties.
\begin{proposition}\label{variety def}
Let $X$ be a proper smooth variety over a perfect field $\k$ and let $\Lambda$ be
a complete local ring with residue field $\k$. Suppose that for any deformation
$\sX \to \Sp A$ of $X$ over an algebra $A \in C_{\Lambda}$, the reduction map
$H^1(\sX,T_{\sX/A}) \to H^1(X,T_{X/\k})$ is surjective. Then the functor of
deformations of $X$ over elements of $C_{\Lambda}$ has the tangent lifting property.
\begin{proof}
The proposition follows directly from (\ref{T1}) and the usual obstruction theory.
\end{proof}
\end{proposition}
\begin{remark}
In the proof of \cite[p.~230]{matsumura86::commut} a simultaneously weaker and
stronger result is obtained. It can be reformulated as saying that if $R$ is a
complete local ring in characteristic zero and the image of $T_R \to
\Hom_{R/\maxid_R}(\maxid_R/\maxid_R^2,R/\maxid_R)$ has dimension $k$, then $R$
is isomorphic to $\pow[S]{t_1,\dots,t_k}$ for some complete subring $S$ of
$R$. (Similar results can be shown in positive and mixed characteristic.)
Applied to the deformation hull $R$ of a smooth and proper variety (still in
characteristic zero) it shows that if $\k$ is the dimension of the image of
$H^1(\sX,T_\sX) \to H^1(X,T_X)$, where $\sX$ is a miniversal deformation, then
$R$ is of the form $\pow[S]{t_1,\dots,t_k}$. Not much seems to be known about
this invariant. Can it, for instance, be trivial for a non-rigid surface $X$?
\end{remark}
The application in characteristic zero of this result (or the corresponding
result for the $T^1$-lifting property) to Calabi-Yau varieties is through two
facts. First that any deformation of a smooth proper variety with $\omega_X$
trivial has a trivial relative $\omega$ and second that the tangent lifting
property is true as it comes down to a lifting property of Hodge cohomology
which is always true. In positive characteristic neither of these facts remain
true. An $\alpha_2$-Enriques surface which has trivial $\omega$ deforms to a
$\Z/2$-Enriques surface whose $\omega$ has order $2$. Similarly Hodge cohomology
does not always have the lifting property (though we do not know of any examples
where the $\omega$ stays trivial in the family). However, if one defines
Calabi-Yau varieties (as does Hirokado, cf.\ \cite{hirokado99::ae+calab+yau}) by
the condition that $\omega_X$ be trivial and $H^i(X,\sO_X)=0$ for $0< i < \dim
X$ then any infinitesimal deformation of a Calabi-Yau variety has trivial
relative $\omega$ in any characteristic. Still it is not clear whether under
that assumption the tangent lifting property is automatic. Hence rather than
formulating a highly conditional result we will give a weaker result that still is
true under rather weak conditions and is the most direct generalisation of the
arguments used in characteristic zero. (For the sake of completeness as well as
comparison we also give the characteristic zero result.) 
\begin{proposition}\label{CY deformations}
Let $X$ be a smooth and proper purely $n$-dimensional variety over a perfect
field $\k$ with $\omega_X$ trivial. Let $F$ be the deformation functor of $X$
over some complete local ring $\Lambda$ with residue field $\k$. When $\ch \k=p
>0$ assume furthermore that $(\Lambda,\maxid_{\Lambda})$ has been given a
divided power structure compatible with the canonical map $\Witt(\k) \to \Lambda$
and the unique divided power structure on $(\Witt,p\Witt)$.

\part If $\ch \k=0$ then $F$ has the tangent lifting property.

\part If $\ch \k=p >0$, $p\Lambda=0$, and $\dim_{\k} H^n_{DR}(X/\k)=\sum_{i=0}^nh^{in-i}_X$ then
$F$ has the divided power tangent lifting property.

\part If $\ch \k=p >0$ and $b_n(X)=\sum_{i=0}^nh^{in-i}_X$ then $F$ has the
divided power tangent lifting property.
\begin{proof}
In all cases we want to prove first that the relative $\omega$ remains trivial
over deformations of $X$ over some rings of $C_{\Lambda}$. If so, then for a
deformation $\sX \to S$ over one of these rings $H^1(\sX,T_{\sX/S})$ equals
$H^1(\sX,\Omega^{n-1}_{\sX/S})$ and then we want to prove that
$H^1(\sX,\Omega^{n-1}_{\sX/S}) \to H^1(X,\Omega^{n-1}_{X/\k})$ is
surjective. Both would follow if we could show that
$H^i(\sX,\Omega^{n-i}_{\sX/S})$ commutes with base change for all $0 \le i \le
n$. For the first part this is \scite[Thm.\ 5.5]{deligne68::theor+lefsc} for all
rings in $C_{\Lambda}$.

The only case where characteristic zero is used in
\lcite{deligne68::theor+lefsc} is for the degeneration of the Hodge to de Rham
spectral sequence and the fact that $H^n_{DR}(-/-)$ commutes with base
change. The only consequence of the degeneration that is used is exactly that
$\dim_{\k} H^n_{DR}(X/\k)=\sum_{i=0}^nh^{in-i}_X$. On the other hand, let us
recall some results from crystalline cohomology: If $A \in C_{\Lambda}$ and $\sX
\to \Sp A$ is a deformation of $X$, and if $(A,\maxid_A)$ has a divided power
structure compatible with that of $\Lambda$ then
$R\Gamma(\sX,\Omega_{\sX/A}^\cdot)$ is isomorphic to
$R\Gamma(X/\Witt)\Tensor_{\Witt}^LA$ resp.\ to
$R\Gamma(X,\Omega_{X/\k}^\cdot)\Tensor_{\k}^LA$ if $p\Lambda=0$ (cf.\
\cite[Cor.~V:3.5.7]{berthelot74::cohom}). Consider first the case when
$p\Lambda=0$. Then using the fact that the $H^*_{DR}(X/\k)$ are flat
$\k$-modules (as all $\k$-modules are) we get that
$H^*_{DR}(\sX/A)=H^*_{DR}(X/\k)\Tensor_{\k}A$ and as also $\dim_{\k}
H^n_{DR}(X/\k)=\sum_{i=0}^nh^{in-i}_X$ Deligne's argument applies. In the
general case, the condition $b_n(X)=\sum_{i=0}^nh^{in-i}_X$ is equivalent to
$H^i_{crys}(X/\Witt)$ being torsion free for $i=n,n+1$ and $\dim_{\k}
H^n_{DR}(X/\k)=\sum_{i=0}^nh^{in-i}_X$. The first part implies that
$H^n_{DR}(\sX/A)=H^n_{crys}(X/\Witt)\Tensor_{\Witt} A$ so again Deligne's argument
applies.
\end{proof}
\end{proposition}
We have now enough results to prove Theorem \ref{Main CY Theorem}.
\begin{proofof}{Theorem \ref{Main CY Theorem}}
This is proved by combining (\ref{CY deformations}), (\ref{TLP}) and specifically Theorem
\ref{Main DP Theorem}
for the case of positive and mixed characteristic deformations.
\end{proofof}
We do not know of any example of a Calabi-Yau variety fulfilling our conditions,
yet having an obstructed deformation space. We shall finish this section by a
global result that shows that such examples would have to be quite isolated.
\begin{proposition}
\part Let $\k$ be a perfect field of positive characteristic $p$, $S$ a
connected $\k$-scheme of finite type and \map fXS a smooth and proper morphism
of pure dimension $n$ that is versal at all of the points of $S$. Assume that
for each $s \in S$ $\dim_{\k(s)}H^n_{DR}(X_s/\k(s))=\sum_{i+j=n}h^{ij}({X_s})$,
that $\dim_{\k(s)}H^n_{DR}(X_s/\k(s))$ is independent of $s$,
and that there is one point $s \in S$ for which $\omega_{X_s}$ is trivial. Then
$S$ is either a smooth $\k$-variety or $S$ is non-reduced at all its points.


\part Suppose instead that $S$ is a $\Witt(\k)$-scheme everywhere versal over $\Witt$
and that $\overline S:=S\Tensor_{\Witt}\k$ is connected. If $b_n(X_s)=\sum _{i+j=n}
h^{i,j}(X_s)$ for all $s\in S$, $\omega_{X_s}$ is trivial for some $s\in S$,
$X_t$ has a formal lifting over a torsion-free base for some $t\in S$ and $S$ is
somewhere smooth, then $S$ is smooth over $\Witt$.
\begin{proof}
For the first part, we begin by noting that the functions $s \mapsto
h^{in-i}_{X_s}$ are constant. For this it is enough to show that they take the
same value under a specialisation but they rise under specialisation under
semi-continuity but by assumption their sum is constant. From this it follows in
particular that sections of $\omega_{X_s}$ lift to a generisation so that the
set of points for which $\omega_{X_s}$ is trivial is open. As it is also closed
it equals $S$ by assumption. Hence we get also that $s \mapsto
\dim_{\k(s)}H^1(X_s,T^1_{X_s})$ is constant on $S$.

Assume now that $S$ is reduced somewhere. Then there is an irreducible
component $Z$ of $S$ meeting the smooth locus $\Reg(S)$ of $S$; since $S$
is connected, it is enough to show that $S$ is smooth at every point $s$
of $Z$. 

Suppose that \map{g}{Y}{D} is miniversal at $d\in D$ with $Y_d\cong X_s$. Then,
after shrinking (in the \'etale topology) $S$ and $D$ appropriately, there is a
smooth surjective morphism $S \to D$ with $s\mapsto d$. Considering the
composite $Z \to S \to D$ shows that $\Reg(D)$ is not empty and that the function
$e\mapsto h^1(Y_e,T^1)$ is constant on $D$. Let $e\in \Reg(D)$.  Then $\dim_e(D)
\ge h^1(Y_e,T^1)$, since $g$ is versal everywhere, and
\begin{displaymath}
\dim_e(D) \le \dim_d(D) \le h^1(Y_d,T^1) = h^1(Y_e,T^1).
\end{displaymath}
So equality holds throughout, and we are done.

As for the last part we begin by noting that $s \mapsto b_n(X_s)$ is locally
constant by smooth and proper base change and hence constant by the
connectedness of $S$. Hence by the first part $\overline S$ is everywhere
smooth. Assume that we can show that for every $s \in \overline S$ a deformation
hull fulfills is tangent smooth (over $\Witt$). Then by Corollary \ref{T-smooth}
the closure of $S\Tensor_{\Witt}\Witt[1/p]$ intersected with $\overline S$ is open and
as it is obviously is closed and non-empty by assumption it equals $S$.

By Proposition \ref{TLP} a deformation hull $R$ is tangent smooth precisely when
the deformation functor fulfills the TLP. To prove it, it suffices, as in the
proof of Proposition \ref{CY deformations}, to show that $H^n_{DR}(\cX/T)$ is
free of rank $b_n(X_s)$ for all deformations $\cX/T$ of $X_s$, $s \in S$. For
this we may work locally and assume that $\overline S$ is the reduction of a
smooth $\Witt$-scheme $S'$ and we may also assume that we have a $\Witt$-homomorphism
$R' := \widehat{\cO}_{S',s} \to R$. By \cite[Prop.\ 3.15]{B-O78} there is a
(unique) divided power structure on $(R',pR')$ compatible with the standard
structure on $\Witt$ (for the structure map $\Witt \to R$). The same thing is true for
$R$ and the map $R' \to R$ is a divided power map. Let us first consider the
crystalline cohomology $R\Gamma(X/S')$. It is a perfect complex of
$\cO_{S'}$-modules (cf.\ \cite[Thm.\ VII:1.1.1]{berthelot74::cohom}). For each
geometric point $\overline{t} \in S$ we then have that
$R^n\Gamma(X/S')\Tensor^L\Witt(\k(\overline{t}))=H^n(X_{\overline{t}/\Witt})$ (cf.\
\cite[Cor.\ V:3.5.7]{berthelot74::cohom}) and as that by assumption is free of
rank $b_n$ we get by Grauert's theorem $R^n\Gamma(X/S')$ is a free
$\cO_{S'}$ and commutes with base change. For any deformation $\cX \to T$ we
have a map $T \to \Spf R$ such that $\cX$ is the pull back of a versal family
and then $H^n(\cX/T)=R^n\Gamma(X/S')\Tensor\cO_T$ again by \cite[Cor.\
V:3.5.7]{berthelot74::cohom}.
\end{proof}
\end{proposition}
\begin{proofof}{Theorem \ref{CY locus}}
The characteristic $p$ part follows directly from the proposition while the
mixed characteristic part follows from it once one has notice that for points in
the closure of the characteristic $0$ locus we have by definition a torsion free
lifting.
\end{proofof}
\end{section}
\bibliography{preamble,abbrevs,alggeom,algebra,ekedahl}
\bibliographystyle{\PRbibstyle}
\end{document}